\documentclass[12pt]{article}
\usepackage{graphicx}
\usepackage{color}
\usepackage{amsfonts}
\usepackage{amsmath}
\usepackage[mathcal]{eucal}
\usepackage{amsthm}
\usepackage{amssymb}
\usepackage{caption}
\usepackage{enumerate}
\usepackage{a4wide}
\usepackage{epstopdf}
\usepackage{mathtools}
\usepackage{soul}
\usepackage[normalem]{ulem}
\usepackage{float}
\usepackage{caption}

\renewcommand{\subset}{\subseteq}

\newtheorem{theorem}{Theorem}

\newtheorem{lemma}[theorem]{Lemma}
\newtheorem{remark}[theorem]{Remark}

\newtheorem{conjecture}[theorem]{Conjecture}

\title{Clique immersions and independence number}

\author{Sebasti\'an Bustamante\thanks{Departamento de Ingenier\'ia Matem\'atica, Universidad de Chile, Chile. SB acknowledges support by ANID Doctoral Fellowship grant 21141116.}, \  Daniel A. Quiroz\thanks{Instituto de Ingenier\'ia Matem\'atica -- CIMFAV, Universidad de Valparaiso, Chile. DAQ acknowledges support from  FONDECYT/ANID Iniciaci\'on en Investigaci\'on Grant 11201251, from   MATH-AMSUD MATH210008 and MATH190013, from FAPESP-ANID Investigaci\'on Conjunta grant 2019/13364-7, and from CMM ANID PIA Basal AFB170001.}, \  Maya Stein\thanks{\,Departamento de Ingenier\'ia Matem\'atica, Universidad de Chile and  Centro de Modelamiento Matem\'atico, Universidad de Chile and CNRS IRL 2807, Chile. MS was supported by ANID Regular Grant 1180830, by MathAmSud MATH190013, by FAPESP-ANID Investigaci\'on Conjunta grant 2019/13364-7, and by ANID Basal AFB170001, ACE210010 and FB210005.}  \ and \ Jos\'e Zamora\thanks{Departamento de Matem\'aticas, Universidad Andr\'es Bello, Chile. JZ acknowledges support  by FONDECYT Regular grant 1180994 and by FAPESP-ANID Investigaci\'on Conjunta grant 2019/13364-7.  } }

\date{}
\begin{document} 
\maketitle

\begin{abstract} 
The analogue of Hadwiger's conjecture for the immersion order
 states that every graph $G$ contains $K_{\chi (G)}$ as an immersion. If true, this would imply that every graph with $n$ vertices and independence number $\alpha$ contains $K_{\lceil \frac n\alpha\rceil}$ as an immersion. 

The best currently known bound for this conjecture is due to
 Gauthier, Le and Wollan, who recently proved that every graph $G$ contains an immersion of a clique on $\bigl\lceil \frac{\chi (G)-4}{3.54}\bigr\rceil$ vertices. 
Their result implies that every $n$-vertex graph with independence number~$\alpha$ contains an immersion of a clique on $\bigl\lceil \frac{n}{3.54\alpha}-1.13\bigr\rceil$  vertices. 

We improve on this result for all $\alpha\ge 3$, by showing that every $n$-vertex graph with  independence number $\alpha\ge 3$ contains an immersion of a clique on~$\bigl\lfloor \frac {n}{2.25 \alpha -f(\alpha)} \bigr\rfloor - 1$ vertices, where $f$ is a nonnegative function.

\end{abstract}

\section{Introduction}

A famous conjecture of  Hadwiger~\cite{H43} states that every graph $G$ of chromatic number $\chi (G)\ge t$ contains  $K_{t}$ as a minor. The conjecture is only known to be true for $t\le 6$, and  probably hard for larger values of $t$, as the proofs for cases $t=5, 6$ already depend on the Four Color Theorem. In general, results of Kostochka~\cite{K84} and, independently, Thomason~\cite{T84} imply that every graph~$G$ with $\chi(G)\ge t$ contains  a clique on $\lceil \frac {t}{f(t)} \rceil$ vertices as a minor,  where $f(t)\in O(\sqrt{\log t})$. Recent breakthroughs of Norin, Postle, and Song~\cite{NPZ20} and of Postle~\cite{P20} have improved on this general bound, but it is still unknown whether $f(t)=c$ for some constant $c$.

It is easy to see that $\chi(G)\ge \lceil\frac n{\alpha(G)}\rceil$ for every $n$-vertex graph $G$, where $\alpha(G)$ denotes the independence number of~$G$. So, 
Hadwiger's conjecture, if true, implies that every $n$-vertex graph contains a clique minor of order at least $\lceil\frac n{\alpha(G)}\rceil$. Duchet and Meyniel~\cite{DM82} conjectured that this holds. Note that the order of this clique minor would be 
 best possible as $G$ could be the disjoint union of cliques. 

Providing evidence for their conjecture, 
Duchet and Mey\-niel~\cite{DM82} proved that every graph~$G$ on $n$ vertices contains a clique minor of order $\lceil\frac{n}{2\alpha(G)-1}\rceil$. There have been several improvements on the order of the clique minor~\cite{BLW11,F10,KPT05,KS07,W07}, the best bound, due to Balogh and Kostochka~\cite{BK11}, being $\lceil\frac{n}{c\alpha (G)}\rceil$, where $c$ is a constant with $c<1.95$. 

The focus of this paper is an analogous result replacing minors  with  immersions. 
A graph~$G$ is said to contain another graph $H$ as an \emph{immersion} if there exists an injective function $\phi\colon V(H)\rightarrow V(G)$ such that:
\begin{enumerate}[(I)]
\item For every $uv\in E(H)$, there is a path in $G$, denoted $P_{uv}$, with endpoints $\phi(u)$ and~$\phi(v)$.
\item The paths in $\{P_{uv} \mid uv\in E(H) \}$ are pairwise edge disjoint.
\end{enumerate}
The vertices in $\phi(V(H))$ are called the \emph{branch vertices} of the immersion. If branch vertices are not allowed to appear as interior vertices on paths $P_{uv}$, the immersion is called {\em strong}.

The minor order and the immersion order are not comparable. The class of planar graphs, while excluding~$K_5$ as a minor, contains all cliques as immersions. On the other hand, the class of graphs with maximum degree at most $d$, while excluding $K_{d+2}$ as an immersion, contains all cliques as minors. 
However, the two relations do share some important similarities. Both of them are well-quasi-orders~\cite{RS04,RS10}, and both notions are weakenings of the topological minor relation. Possibly inspired by such similarities, Lescure and Meyniel~\cite{LM89}  proposed an analogue of Hadwiger's conjecture for strong immersions. Later, Abu-Khzam and Langston~\cite{AKL03} weakened their conjecture to the following form.

\begin{conjecture}[\cite{AKL03, LM89}] \label{hadimm}
Every graph $G$ with $\chi(G)\ge t$ contains $K_{t}$ as an immersion.
\end{conjecture}

This conjecture and its strong version have received much attention recently, and have been tackled with more success than their minor counterpart. The cases $1\le t \le 4$  follow from the fact that Haj\'os' Conjecture 
 is true for these cases~\cite{D52}. The cases $5\le t \le 7$ of Conjecture~\ref{hadimm} were established by Lescure and Meyniel~\cite{LM89} and by DeVos, Kawarabayashi, Mohar, and Okamura~\cite{Detal10}. 

For general values of $t$, the first linear lower bound for Conjecture~\ref{hadimm} was given by DeVos, Dvo\v{r}\'ak, Fox, McDonald, Mohar, and Scheide~\cite{Detal14}. They proved that every graph~$G$ contains an immersion of a clique on $\bigl\lceil \frac{\chi (G)}{200} \bigr\rceil$ vertices. Dvo\v{r}\'ak and Yepremyan~\cite{DY17} improved this bound to $\bigl\lceil \frac{\chi (G)-7}{11} \bigr\rceil$. The best currently known lower bound for Conjecture~\ref{hadimm} is due to Gauthier, Le and Wollan~\cite{GLW17} who showed that every graph $G$ contains an immersion of a clique on $\bigl\lceil \frac{\chi (G)-4}{3.54} \bigr\rceil$.  This implies the following.

\begin{theorem}[Gauthier, Le and Wollan~\cite{GLW17}]\label{wow}
Every $n$-vertex graph $G$ contains  an immersion of a clique on  $\bigl\lceil \frac{n}{3.54\alpha(G)}-1.13\bigr\rceil$ vertices.
\end{theorem}

It seems natural to ask whether the bound from Theorem~\ref{wow} can be improved without necessarily improving Gauthier, Le and Wollan's underlying result that relates the chromatic number and the size of clique immersions,  as Duchet and Meyniel did in the context of graph minors. 
The first attempt in this direction (actually earlier than~\cite{GLW17}) has been carried out by Vergara~\cite{V17}. She conjectured that 
 every $n$-vertex graph with  independence number 2 contains an immersion of $K_{\lceil \frac n2\rceil}$, and showed that this conjecture is equivalent to Conjecture~\ref{hadimm} for graphs of independence number 2.
 In support of her conjecture, Vergara  proved that every graph on $n$ vertices and independence number 2  contains a $K_{\lceil \frac n3\rceil}$-immersion. Her result was improved by Gauthier, Le and Wollan~\cite{GLW17} as follows. (See~\cite{Q20} for other results on this conjecture.)

\begin{theorem}[Gauthier, Le and Wollan~\cite{GLW17}]\label{dos}
Every $n$-vertex  graph with independence number~$2$ contains $K_{2\lfloor \frac n5 \rfloor}$ as an immersion.
\end{theorem}

Extending the conjecture of Vergara, we propose the following conjecture for graphs of arbitrary independence number.

\begin{conjecture}\label{ourconj}
Every $n$-vertex  graph $G$ contains $K_{\lceil \frac n{\alpha(G)}\rceil}$ as an immersion.
\end{conjecture}

This conjecture is best possible, since (as mentioned earlier for the case of minors) the graph  $G$ could be the disjoint union of cliques. 

Our main result improves on Theorem~\ref{wow} for every $G$ with $\alpha(G)\ge 3$, giving more evidence in support of Conjecture~\ref{ourconj}, and thus of Conjecture~\ref{hadimm}.

\begin{theorem}\label{all}
Let $G$ a graph on $n$ vertices  of independence number 
$\alpha \geq 3$.  Then $G$ contains an immersion of a clique on at least $\lfloor  \frac {n}{2.25\alpha  -f(\alpha)} \rfloor - 1$ vertices, where $f(3)=f(4)=2.25$ and $f(\alpha)=2.25+\sum_{i=6}^{\alpha+1} \frac 1{2i}$, if $\alpha\ge 5$.
\end{theorem}

In particular, Theorem~\ref{all}  implies that any $n$-vertex graph with independence number $\alpha\le 10$ contains an immersion of a clique on $\lfloor \frac {n}{2 \alpha} \rfloor - 1$ vertices.

In the literature, several other types of immersions have been distinguished (apart from immersions and strong immersions). If all the paths~$P_{uv}$ from~(I) have length at most $k$, then~$H$ is called a \emph{$k$-immersion} of $G$. If  all the paths $P_{uv}$ have odd length then~$H$ is an \emph{odd immersion}. 

Actually, Theorems~\ref{wow} and~\ref{dos} as well as   the precursor results give  strong immersions. Moreover, one can read from
the proof of Theorem~\ref{dos}  that the immersion from that theorem is a strong odd $5$-immersion. While no such information can be deduced for Theorem~\ref{wow}, it is not difficult to verify that the proof of our Theorem~\ref{all} does imply a stronger statement. 

\begin{remark}\label{more}
The immersion from Theorem~\ref{all} is  a strong odd $(2\alpha -1)$-immersion.
\end{remark}

In particular, this means that our result not only gives further evidence for Conjectures~\ref{hadimm} and~\ref{ourconj}, but also for a stronger conjecture of Churchley~\cite{C17} stating that every graph of chromatic number at least~$t$  contains $K_{t}$ as an odd immersion.

The structure of the paper is as follows. In the next section (Section~\ref{secaalpha}) we introduce the auxiliary graphs $A_\alpha$ and prove a crucial lemma, Lemma~\ref{parti}. This lemma will be used at the end of the proof of Theorem~\ref{all}, which will be given  in Section~\ref{secall}. 

An extended abstract~\cite{nos} announcing a weaker version of Theorem~\ref{all}  includes an alternative proof for the case $\alpha(G)=3$.

\section{The graph $A_\alpha$ and its minimal cuts}\label{secaalpha}

In this section we give key definitions and prove some crucial ingredients for the proof of Theorem~\ref{all}.

\subsection{Minimal cuts and blow-ups}

We first need to introduce some standard notation on (vertex-)cuts. Let $G$ be a graph. For $A,B\subseteq V(G)$, an \emph{$A$-$B$ cut} is a set $C\subseteq V(G)$ such that $G\setminus C$ contains no path from $A\setminus C$ to $B\setminus C$. (Slightly abusing notation, we might sometimes write $a$-$B$ cut if $A=\{a\}$.) 
If $C$ is an $A$-$B$ cut we also say $C$ is \emph{a cut separating $A$ and $B$}. We say $C$ is \emph{minimal} if for every $x\in C$, we have that $C\setminus\{ x\}$ is not an $A$-$B$ cut. 

Let $\mathbb{N}$ denote the set of nonnegative integers. For a graph $G$ and $f\colon V(G)\rightarrow\mathbb{N} $ a function, an \emph{$f$-blow-up}  of $G$ is a graph  $B$ that can be obtained from $G$ by replacing each vertex $v\in V(G)$ with an independent set $B(v)$ satisfying $|B(v)|=f(v)$, and $B(v)\cap B(u)=\varnothing$ if $v\ne u$,  and having the edge $xy$, for $x\in B(v), y\in B(u)$, if and only if $uv\in E(G)$.  We  let $B(X)=\bigcup_{v\in X}B(v)$ for every $X\subseteq V(G)$.

We will need  the following lemma. 
\begin{lemma}\label{blowcut}
Let $G$ be a graph, let $f\colon V(G)\rightarrow \mathbb{N}$, and let $X, Y\subseteq V(G)$. Let $B$ be an $f$-blow-up of~$G$, and let $C_B$ a minimal $B(Y)$-$B(X)$ cut in $B$.  Then there  is an $X$-$Y$ cut $C$ in $G$ such that  $C_B=\bigcup_{v\in C}B(v)$.
\end{lemma}
\begin{proof}
Choose $Z\subseteq \{v\in V(G) :f(v)=0\}$ minimal such that 
$C=\{v\in V(G) : B(v)\cap C_B\neq\varnothing\}\cup Z$ is  an $X$-$Y$ cut in $G$. Then clearly $C_B=\bigcup_{v\in C}B(v)$.
\end{proof}

\subsection{The  auxiliary graph $A_\alpha$}

We now define a family of graphs $\{ A_\alpha \}_{\alpha\ge 2}$. These graphs will represent, in a simplified way, the structure that arises when trying to find the desired immersion for our theorem by induction. 


There is another auxiliary graph involved in our proof, and this is a blow-up of~$A_\alpha$. The function defining the blow-up varies with the graph $G$ of independence number $\alpha$ from Theorem~\ref{all}. So it is good to think of each vertex (except $d_\alpha$) of the auxiliary graph $A_\alpha$ as representing a subset of $V(G)$, and each edge representing a complete bipartite subgraph. Yet, the blow-up of $A_\alpha$ is not necessarily isomorphic to a subgraph of $G$. The exact relation between the blow-up of $A_\alpha$ and $G$ will be clarified later. For the final structural analysis however, it is easier (and sufficient) to restrict our attention to~$A_\alpha$.

For a given $\alpha\ge 2$ we define  $A_\alpha$ as the graph having the following vertices:
\begin{itemize}\itemsep -2pt
\item a vertex $d_\alpha$;

\item  a vertex $d_S$ for every nonempty $S\subseteq \{1,\dots ,\alpha-1\}$; and

\item a vertex $x_S$ for every $S \subseteq \{1,\dots ,\alpha\}$ with $|S|\ge 2$;
\end{itemize}

\noindent
and the following  edges:
\begin{itemize}\itemsep -2pt
\item the edge $d_\alpha x_S$, for every $S \subseteq \{1,\dots ,\alpha\}$ with $\alpha\in S$,  and $|S|\ge 2$;

\item the edge $d_Sx_T$, for each pair $S,T\subseteq \{1,\dots ,\alpha\}$ with $\alpha\notin S$, $|T|\ge 2$, and $S\cap T\ne\varnothing$;

\item the edge $x_Sx_T$, for each pair $S,T\subseteq \{1,\dots ,\alpha\}$ with $|S|, |T|\ge 2$, and $S\cap T\ne\varnothing$.
\end{itemize}

We write $$D_\alpha:=\{d_S \mid S \subseteq \{1,\dots ,\alpha-1\}, S\ne \varnothing\},$$ and observe that 
\begin{equation}\label{neighbs}
N_{A_\alpha}(d_T)\subseteq N_{A_\alpha}(d_S)
\end{equation}
for all $T\subseteq S$ with $d_T, d_S \in D_\alpha$.

Since we plan to use Menger's theorem later for finding the connections between the new branch vertex and the old ones, we will be interested in minimal cuts of the original graph~$G$. Here we investigate the minimal cuts of the graphs $A_\alpha$, and later use Lemma~\ref{blowcut} to translate what we find here into information about the cuts of $G$. The following lemma contains all the important properties of minimal cuts we will use. 

\begin{lemma}\label{Aalpha}
Let $\alpha\ge 2$ and let $C$ be a minimal $d_\alpha$-$D_\alpha$ cut in $A_\alpha$ such that $C\ne \{d_\alpha\}$. Then all of the following hold, for each $S\subset \{1,\dots ,\alpha-1\}$.
\begin{enumerate}[(a)]
\item\label{aalphad}
If  $d_S\notin C$, then $d_T\notin C$ for every nonempty $T\subset S$.
\item\label{aalphax}
If  $x_{S\cup\{\alpha\}}\notin C$, then $x_{T\cup\{\alpha\}}\notin C$ for every nonempty $T\subset S$.
\item\label{aalpha1}
If $d_S\notin C$, and if $x_{T\cup\{\alpha\}}\notin C$ for some nonempty $T\subseteq \{1,\dots ,\alpha-1\}$, then $x_R \in C$ for every $R\subseteq \{1,\dots ,\alpha-1\}$ with $R\cap S\ne\varnothing\ne R\cap T$ and $|R|\ge 2$.
\item
\label{aalpha2}
We have $d_S\notin C$ if and only if $x_{T\cup\{\alpha\}}\in C$ for every nonempty $T\subseteq \{1,\dots ,\alpha-1\}$ with $T\cap S\ne \varnothing$.
\end{enumerate}
\end{lemma}

\begin{proof}
We first prove~\eqref{aalphad}. Assume there is $T\subset S$ with $d_T\in C$. Since $C$ is minimal, $A_\alpha\setminus (C-d_T)$ contains a $d_T$-$d_\alpha$ path. But then, by~\eqref{neighbs}, either  $A_\alpha\setminus C$ contains a $d_S$-$d_\alpha$ path, which is impossible, or $d_S\in C$, which is as desired. This proves~\eqref{aalphad}.

Property~\eqref{aalphax} can be proved using similar arguments as for~\eqref{aalphad} since $N_{A_\alpha}(x_{T\cup\{\alpha\}})\subseteq N_{A_\alpha}(x_{S\cup\{\alpha\}})$. And property~\eqref{aalpha1} holds since $d_Sx_Rx_{T\cup \{\alpha\}}d_\alpha$ is a path.

Let us now prove~\eqref{aalpha2}.
Suppose first that we have $d_S\notin C$. If $T\cap S\ne \varnothing$, consider the path $d_Sx_{T\cup\{\alpha\}}d_\alpha$. This path has to meet $C$, and thus we have  $x_{T\cup\{\alpha\}}\in C$.

Suppose now that $x_{T\cup\{\alpha\}}\in C$ for every nonempty $T$ with $T\cap S\ne \varnothing$. For a contradiction, assume $d_S\in C$. Since~$C$ is minimal, there must be a $d_S$-$d_\alpha$ path in $A_\alpha\setminus(C-d_S)$. Let $P=d_Sx_{R_k}\dots x_{R_1}x_{R_0}d_\alpha$ be a minimum length path among these. By hypothesis we have that $k\ge 1$. 

Let us show that there exists $U\neq\varnothing$ such that $U\subseteq R_k$ and $x_{U\cup\{\alpha\}}\in C$. If $k=1$, then by definition of $P$, we have $S\cap R_1\ne \varnothing$. So we can take $U=S\cap R_1$ (notice that $x_{U\cup \{\alpha\}}\in C$ by hypothesis). If $k\ge 2$, we note that since $P$ is of minimum length, the path $d_Sx_{R_k}x_{R_k\cup\{\alpha\}}d_\alpha$ meets $C$, which implies $x_{R_k\cup\{\alpha\}}\in C$. So we can take $U=R_k$.

By minimality of $C$, the graph $A_\alpha\setminus(C-x_{U\cup\{\alpha\}})$ contains a $D_\alpha$-$d_\alpha$ path that passes through $x_{U\cup\{\alpha\}}$. Let $x_W$, or alternatively $d_W$, be the vertex before $x_{U\cup\{\alpha\}}$ on this path, and $P'$ the subpath of this path that joins $D_\alpha$ with either $x_W$ or $d_W$. If we had $\alpha\in W$, the vertex before $x_{U\cup\{\alpha\}}$ would not be $d_W$, but rather $x_W$ and adjacent to $d_\alpha$, and there would be a $D_\alpha$-$d_\alpha$ path that avoids $C$, a contradiction. Then we have $\alpha\notin W$ and, since $U\subseteq R_k$, we have  $\varnothing \ne W\cap U\subseteq W\cap R_k$. Then  $P'x_{R_k}\dots x_{R_1}x_{R_0}d_\alpha$ is or contains a $D_\alpha$-$d_\alpha$ path that avoids $C$, a contradiction.
\end{proof}

We are now ready to show the crucial lemma for the proof of Theorem~\ref{all}.

\begin{lemma}\label{parti}
Let $\alpha\ge 2$ and let $C$ be a minimal $d_\alpha$-$D_\alpha$ cut in $A_\alpha$ such that $C\ne \{d_\alpha\}$. Then there are two (possibly empty) disjoint sets $J_1$ and $J_2$, such that $J_1\cup J_2=\{1,\dots ,\alpha -1\}$ and
\begin{enumerate}[(i)]
\item $x_S\in C$ for every $S\subseteq \{1,\dots ,\alpha -1\}$ with $J_1\cap S\neq\varnothing\neq J_2\cap S$;
\item $x_{S\cup\{\alpha\}}\in C$ for every $S\subseteq \{1,\dots ,\alpha -1\}$ with $J_2\cap S\neq\varnothing$;
\item\label{vacio} If $J_2=\varnothing$ then $d_S\in C$ for every nonempty $S\subseteq \{1,\dots ,\alpha -1\}$.
\end{enumerate}
\end{lemma}

\begin{proof}
Let $J_1$ be the smallest subset of $\{1,\dots ,\alpha -1\}$ such that  for every nonempty $T \subseteq \{1,\dots ,\alpha -1\}$ with $x_{T\cup \{\alpha\}}\notin C$ we have $T\subseteq J_1$. Let $J_2$ be the smallest subset of $\{1,\dots ,\alpha -1\}$ such that  for every $S\subseteq \{1,\dots , \alpha -1\}$ with $d_S \notin C$ we have $S\subseteq J_2$. Clearly \emph{(iii)} holds for this choice.

By Lemma~\ref{Aalpha}~\eqref{aalpha2}, we know that if $d_S \notin C$ then $x_{T\cup \{\alpha\}}\in C$ for every $S,T\subseteq \{1,\dots ,\alpha -1\}$ with $S\cap T\ne \varnothing$. By the choice of $J_1$ and $J_2$ as minimal, this implies that $J_1\cap J_2=\varnothing$. 
Moreover, if there exists $i\in \{1,\dots ,\alpha -1\}\setminus (J_1\cup J_2)$, then we have  $x_{\{ i\}\cup\{\alpha\}}\in C$ and $d_{\{ i\}}\in C$, which contradicts Lemma~\ref{Aalpha}~\eqref{aalpha2}. Thus, we also have  $J_1\cup J_2=\{1,\dots ,\alpha -1\}$.

By Lemma~\ref{Aalpha}~\eqref{aalphad} and~\eqref{aalphax}, it follows that for every $i\in J_1$ and every $j\in J_2$ we have $x_{\{i\}\cup\{\alpha\}}\notin C$ and $d_{\{ j\}}\notin C$. This implies, by Lemma~\ref{Aalpha}~\eqref{aalpha1}, that for every $i\in J_1$ and every $j\in J_2$ we have $x_S\in C$ for every $S\subseteq \{1,\dots ,\alpha -1\}$ with $i,j\in S$. Thus \emph{(i)} holds. Moreover,   Lemma~\ref{Aalpha}~\eqref{aalpha2} implies that for every $j\in J_2$ we have $x_{S\cup\{\alpha\}}\in C$ for every $S\subseteq \{1,\dots ,\alpha -1\}$ with $j\in S$, and therefore, also \emph{(ii)} holds. 
\end{proof}

\section{The proof of Theorem~\ref{all}}\label{secall}

We first need the following elementary lemma.

\begin{lemma}\label{elem}
Let $\beta_1=1$, $\beta_2=\frac{5}2$, $\beta_3=\frac 92$, $\beta_4=\frac {27}4$ and $\beta_i=\frac94 (i-1) - \sum_{j=6} ^{i+1} \frac 1{2j}$  for $i\ge 5$. Let~$s,t$ and $i\ge 4$ be integers with $s+t=i$, and $2\le s\le t\le i-2$. Then we have $$\beta_i\ge  \beta_s+\beta_t +\frac{2i -1}i.$$
\end{lemma}
\begin{proof}
If $i=4$, then $s=t=2$ and the result holds with equality. So we assume $i\ge 5$.  If $s=2$, notice that we have $$\beta_{i}=\beta_{i-2}+\frac92 -\frac 1{2i}- \frac 1{2(i +1) } > \beta_2 +\beta_{i-2} +\frac{2i -1}i ,$$ as desired. When $s=3$ we must have $i\ge 6$, and thus $$\beta_{i}\ge \beta_{i-3}+\frac{27}4 - \frac 1{10} -\frac 1{2i } - \frac 1{2(i +1) } > \beta_3 +\beta_{i-3} +\frac{2i -1}i .$$The case $s=4$ follows analogously. When $s\ge 5$, it is not hard to see that we have
\begin{equation*}
\beta_{i} = \beta_{t}+\frac {9s}4- \sum_{j=t+2}^{i+1}\frac 1{2(i +1)} >  \beta_t+  \beta_s + \frac{2i -1}i ,
\end{equation*}
which is as desired.
\end{proof}

Theorem~\ref{all} follows immediately from Lemma~\ref{elem} together with the following result. 
\
 \begin{theorem}\label{thmX}
 Let  $\beta_1,\beta_2,\beta_3, \ldots$ be rationals such that $\beta_1=1,\beta_2=\frac 52$, while for all $i\ge 3$ \\ 
\phantom{yap} $\bullet$ \, $\beta_i\ge \beta_{i-1}+2$ and \\
\phantom{yap} $\bullet$ \,  for all  $s, t\in\{2, \ldots, i-2\}$ with $s+t=i$  we have $\beta_i\ge \beta_s+\beta_t +\frac{2i-1}i.$ \\
Let $G$ be a graph on~$n$ vertices of independence number at most $\alpha$. Then $G$ contains an immersion of a clique on at least $\lfloor \frac n{\beta_\alpha} \rfloor - 1$ vertices. 
\end{theorem}




\begin{proof}
We proceed by  induction on $\alpha$. The statement is trivially true for $\alpha=1$, and true for $\alpha=2$ by Theorem~\ref{dos}. We now prove it for~$\alpha\ge 3$, using induction on $n$. It is easy to see that the statement holds for all $1\le n<3\beta_\alpha$ as we then only need to find a $K_1$-immersion.

Note that we can assume that $\alpha (G)=\alpha  <n$, and consider a maximum independent set $I=\{a_1, \dots, a_\alpha\}$ of~$G$. By  induction on $n$, 
we know that $G-I$ contains an immersion  of a clique on $\lfloor \frac{n-\alpha}{\beta_{\alpha(G-I)}} \rfloor - 1$  vertices. In particular, since  $\beta_{\alpha(G-I)}\le \beta_{\alpha(G)}$, $G-I$ has an immersion  of a clique on
$\lfloor \frac{n-\alpha}{\beta_{\alpha}} \rfloor - 1$ vertices, with set of branch vertices~$M$. 
Our plan is to  add a new branch vertex to this immersion, or, failing this, to find  some other immersion of  a clique on at least $\lfloor \frac n{\beta_\alpha} \rfloor - 1$ vertices in $G$. If we succeed in adding a new branch vertex to the immersion, we will have constructed an immersion on $\lfloor \frac{n-\alpha}{\beta_{\alpha}} \rfloor - 1+1\ge \lfloor \frac{n}{\beta_{\alpha}} \rfloor - 1$ vertices, which is as desired. 

Set $Q:=V(G)\setminus (I\cup M)$ and note that
\begin{equation}\label{Qall}
|Q|\ = \ n-\alpha -(\Big\lfloor \frac{n-\alpha}{\beta_\alpha} \Big\rfloor - 1)\ > \ \frac{\beta_\alpha-1}{\beta_\alpha}n - \alpha.
\end{equation}

Set $N_i:=N(a_i)$ and $\bar N_i:= V(G)\setminus (N(a_i)\cup \{a_i\} )$ for $i=1,\dots ,\alpha$. Since $G$ has independence number~$\alpha$, we know that if $S\subseteq\{ 1,\dots,\alpha\}$ then $G[\bigcap_{i\in S}\bar N_i]$ has independence number at most $\alpha-|S|$. Hence, we may assume that 
\begin{equation}\label{nonneighS}
\Bigl| \bigcap_{i\in S} \bar N_i\Bigr| <\frac{\beta_{\alpha-|S|}}{\beta_\alpha}\, n, \mbox{ for every } S\subseteq\{ 1,\dots,\alpha\}, \, 1\le|S|\le \alpha -1
\end{equation}
as otherwise, by induction on $\alpha$, there would be some $S$ such that $G[\bigcap_{i\in S} \bar N_i]$ contains an immersion of $K_{\lfloor \frac{n}{\beta_{\alpha}} \rfloor -1}$.

Let us prove next that we have 
\begin{equation}\label{neeew}
\text{$|  M\cap \bar N_i | < |  Q\cap N_i \cap \bigcup_{j \neq i} N_j|$ for each $i \in \{1,\dots ,\alpha \}$.} 
\end{equation}
Indeed for a contradiction, suppose there exists $i$ such that $| M\cap \bar N_i  | \geq |  Q\cap N_i \cap \bigcup_{j \neq i} N_j |$. Notice that $Q = (\bigcap_{j \neq i} Q\cap \bar N_j) \cup (  Q\cap N_i \cap  \bigcup_{j\neq i} N_j ) \cup ( Q\cap \bar N_i )$. Since $\bigcap_{j \neq i}\bar  N_j$ is a clique, we obtain
\begin{equation*}
\begin{split}
|Q| &\leq \Big\lfloor \frac{n}{\beta_\alpha} \Big\rfloor -2 +| M\cap \bar N_i  |+ |Q\cap \bar N_i |\\
 &\leq \frac{n}{\beta_\alpha}  -2 + |\bar N_i\setminus I| \\ 
&< \frac{n}{\beta_\alpha} -\alpha -1 + \frac{\beta_{\alpha -1}}{\beta_\alpha}n \\
&= \frac{\beta_{\alpha-1}+1}{\beta_\alpha}n -\alpha -1 ,  
\end{split}
\end{equation*}
where the third inequality comes from~\eqref{nonneighS}. This contradicts \eqref{Qall}, since $\alpha\ge 3$ and $\beta_\alpha\ge \beta_{\alpha-1}+2$.

Without loss of generality we can assume $| M\cap \bar N_\alpha | =\min_{ i\in \{ 1,\dots ,\alpha \} } |M\cap \bar N_i|$ and so, by~\eqref{neeew}, 
\begin{equation}\label{boundapares}
| M\cap \bar N_\alpha  | < \Bigl|   Q\cap N_i \cap  \bigcup_{j \neq i} N_j\Bigr|, \mbox{ for every } 1\le i\le \alpha.
\end{equation}

This assumption also implies
\begin{equation}\label{pendiente}
|M\cap N_\alpha | \ge \frac{1}{\alpha}|M|,
\end{equation}
since every vertex of $M$ is a neighbor of some vertex in $I$ (because $I$ is a maximum independent set).

 We  claim that if 
 $I_1$ and $I_2$  are disjoint sets with $ I_2 \ne \varnothing$  and $I_1\cup I_2=\{1,\dots, \alpha-1\}$, then 
 \begin{equation}\label{claimst}
 \Bigl| Q\cap \Bigl(\bigcup_{i\in I_1\cup \{\alpha\}} N_i\Bigr)\cap \Bigl(\bigcup_{j\in I_2} N_j\Bigr)\Bigr|\ge |M\cap \bar N_\alpha  |. 
 \end{equation}
For the case $I_1=\varnothing$ and for the case $|I_2|=1$, \eqref{claimst} follows from \eqref{boundapares}. So we assume $I_1\ne \varnothing$ and $|I_2|\ge2$, and note that we have
\begin{equation*}
\begin{split}
\Bigl| Q\cap \Bigl(\bigcup_{i\in I_1\cup \{\alpha\}} N_i\Bigr)\cap \Bigl(\bigcup_{j\in I_2} N_j\Bigr) \Bigr| &\ge n - |M|- \Bigl|\bigcap_{i\in I_1\cup \{\alpha\}}\bar N_i\Bigr|- \Bigl|\bigcap_{j\in I_2}\bar N_j\Bigr|  \\
&> \Bigl(\beta_\alpha -1 - \beta_{\alpha -(|I_1|+1)} -\beta_{\alpha-|I_2|}\Bigr) \frac n{\beta_\alpha}\\
&\ge \frac{\alpha-1}\alpha\cdot\frac n{\beta_\alpha} \\ & > |M\cap \bar N_\alpha |,
\end{split}
\end{equation*}
where the second inequality comes from \eqref{nonneighS} and the fact that $|M|< \frac n{\beta_\alpha}$; the third inequality comes from the properties of the $\beta_i$; and the last inequality comes from~\eqref{pendiente} and (again) the fact that $|M|< \frac n{\beta_\alpha}$. Thus we have proved~\eqref{claimst}. 

Let $\mathcal{P}_G$ be the set of all paths in $G$ starting in $a_\alpha$ and ending in $M\cap \bar N_\alpha $ which alternate between the set $M\cup Q$ and the set $I=\{a_1,\dots ,a_\alpha\}$ (i.e., each edge in such a path is incident to one vertex from each set) and have all internal vertices in $Q\cup I\setminus\{a_\alpha\}$\footnote{Note that the definition of $\mathcal P_G$ immediately implies that if $P \in \mathcal{P}_G$ then the length of $P$ is odd and is at most $2\alpha-1$.}. Our goal is to show that there exists a set of $|M \cap\bar N_\alpha |$ edge-disjoint paths in $\mathcal{P}_G$,  each having a different endpoint in $M\cap \bar N_\alpha$. Indeed, if we find such a set, then,  as each path from $\mathcal{P}_G$ is internally disjoint from $M$, and edge-disjoint from the paths occupied by the immersion with branch vertices $M$,  the vertices in $M\cup\{a_\alpha\}$ are  the branch vertices of an immersion of the desired size, which finishes the proof.

We  define the following sets and the function $f\colon V(A_\alpha)\rightarrow \mathbb{N}$ (where $A_\alpha$ is the graph defined in Section~\ref{secaalpha})  as follows:
\begin{itemize}\itemsep -2pt

\item  $M_S = (M\cap \bigcap_{i\in S}N_i)\setminus \bigcup_{j\in \{1,\dots ,\alpha\}\setminus S}N_j$ and $f(d_S) = |M_S|$, for every $S\subseteq \{1,\dots ,\alpha-1\}$;

\item $Q_S=(Q\cap \bigcap_{i\in S}N_i)\setminus \bigcup_{j\in \{1,\dots ,\alpha -1\}\setminus S}N_j$ and $f(x_S) = |Q_S|$, for every $S \subseteq \{1,\dots ,\alpha\}$ with $|S|\ge 2$;

\item $f(d_\alpha)= n$.
\end{itemize}





 Let $B$ be an $f$-blow-up of $A_\alpha$ and for 
$v\in V(A_\alpha)\setminus \{d_\alpha\}$
let 
$\pi_v:  B(v) \to D$
a bijection where 
$D=M_S$ if $v=d_S$ and $D=Q_S$ if $v=x_S$
for some $S$.  Let $\mathcal{P}_B$ be the set of all paths in $B$ with one endpoint in~$B(d_\alpha)$ and the other in $\bigcup_{\emptyset\neq S\subseteq\{1,\dots ,\alpha -1\}}B(d_S)$, and internally disjoint from these two sets.

We claim (and will prove below) that 
\begin{equation}\label{set_of_PB_paths}
\text{there is  a set $\mathcal R\subseteq \mathcal{P}_B$ of $|M \cap\bar N_\alpha |$ vertex-disjoint paths in $B$.}
\end{equation}
If this is true, then we can also find
 a subset of $\mathcal{P}_G$ containing $|M \cap\bar N_\alpha |$ edge-disjoint paths. Indeed, let $R=v_0v_1\ldots v_t \in\mathcal R$
with $v_0 \in B(d_\alpha)$, $v_t \in B(d_{S_t})$ for some non-empty $S_t \subseteq \{1,2,\ldots , \alpha-1\}$, and, for $1\le i<t$, $v_i \in B(x_{S_i})$ for some $S_i \subseteq \{1,2, \ldots , \alpha\}$ with $|S_i|\ge 2$. Notice that we can take all the paths in $\mathcal{R}$ (and in particular $R$) to be induced.
To the path $R$ we assign a path $P_R\in \mathcal{P}_G$, which is defined as follows:
$$ P_R= a_\alpha  \pi_{x_{S_1}}(v_1) a_{i_1} \ldots \pi_{x_{S_{t-1}}}(v_{t-1}) a_{i_{t-1}} \pi_{d_{S_t}}(v_t),$$
where for $1\le j< t$, $i_j$ is the minimum element of the set $S_j \cap S_{j+1}$. Since we took $R$ to be induced, $S_i\cap S_{j+1}=\varnothing$ for $i<j$, and so $P_R$ is indeed a path.
See Figure~\ref{theone} for an illustration.




\begin{figure}[h]
 \centering
 \captionsetup{justification=centering}
 \bigskip
 \includegraphics[height=3.5 in]{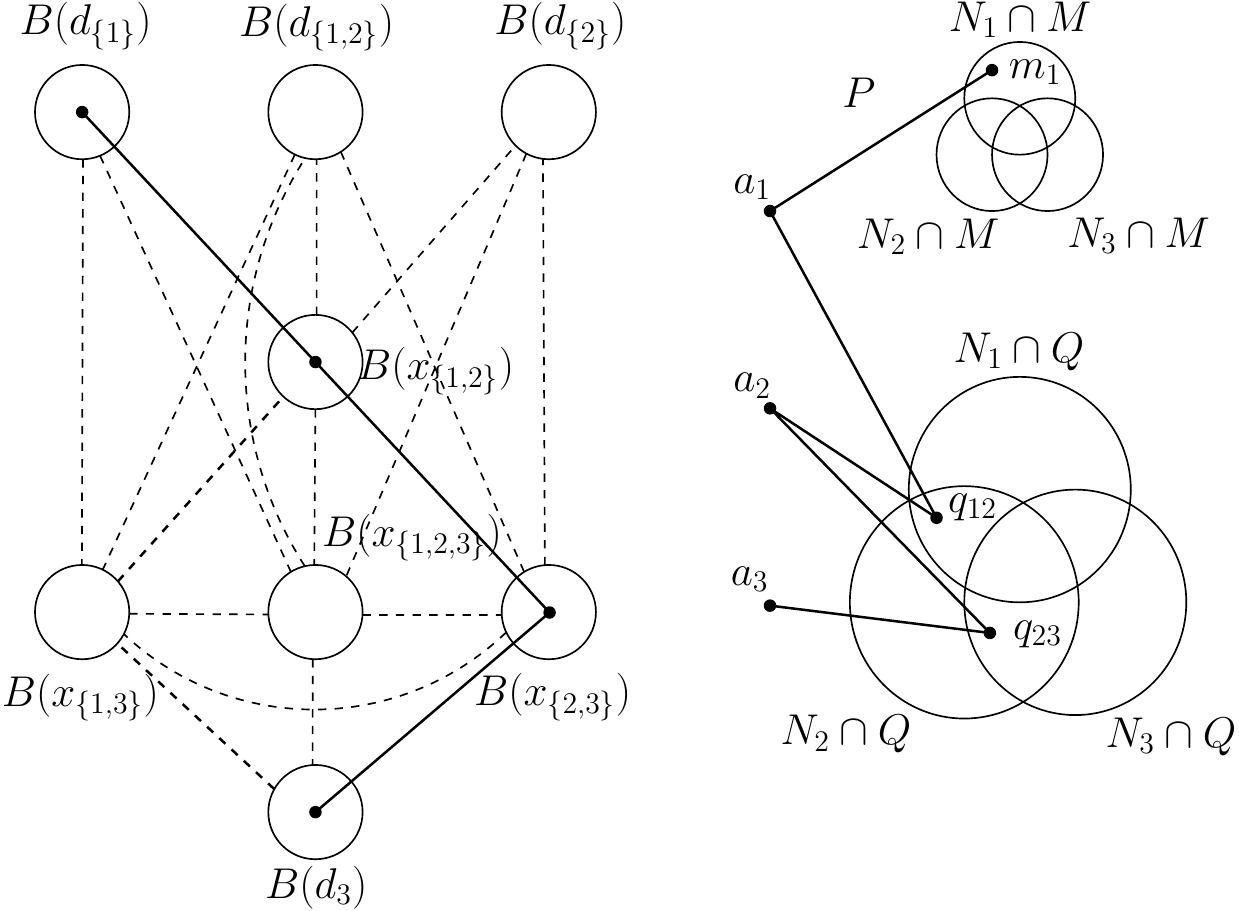}
 \medskip
 \caption{
 A path $R\in \mathcal R$ going through  sets $B(d_3)$, $B(x_{\{ 2, 3\}})$, $B(x_{\{1,2\}})$, $B(d_{\{1\}})$ is assigned a path $P_R=a_3 q_{23} a_2 q_{12} a_1 m_1$ with $q_{23}\in Q_{\{2,3\}}
$, $q_{12}\in Q_{\{1,2\}}
$ and $m_1\in M_{\{1\}}
$.
Dashed edges show 
possible adjacencies between sets $B(v)$ in the graph $B$.
}
  \label{theone}
\end{figure}

 Since the paths in $\mathcal R$ are vertex-disjoint, and the functions $\pi_v$ are bijections, the paths in $\{P_R : R\in\mathcal R\}$ can only intersect in vertices of $I$. That is, they are pairwise edge-disjoint, which is as desired.

 So it only remains to prove~\eqref{set_of_PB_paths}. For this, note that by
Menger's Theorem~\cite{M27}, it suffices to show that every 
minimal cut $C_B$ separating  $B(d_\alpha)$ from $\bigcup_{\emptyset\neq S\subseteq\{1,\dots , \alpha -1 \}}B(d_S)$ fulfills $|C_B|\ge |M\cap \bar N_\alpha|$. Lemma~\ref{blowcut} implies that 
for any such $C_B$,  there  is a  $d_\alpha$-$D_\alpha$ cut $\hat C$ in $A_\alpha$ such that $C_B=\bigcup_{v\in \hat C}B(v)$.
Clearly, we will have $|\bigcup_{v\in \hat C}B(v)|\ge |\bigcup_{v\in C}B(v)|$ for some minimal $d_\alpha$-$D_\alpha$ cut  $C$ in $A_\alpha$. This, and our choice of $f(d_\alpha)$, tells us that we can finish the proof by showing that 

\begin{equation}\label{fini}
\text{$|\bigcup_{v\in C}B(v)|\ge |M\cap \bar N_\alpha|$ for every minimal $d_\alpha$-$D_\alpha$ cut $C$ in $A_\alpha$ with $C\ne \{d_\alpha\}$.}
\end{equation}

In order to see~\eqref{fini}, we take any such cut $C$ and consider the sets $J_1$ and $J_2$ given by Lemma~\ref{parti}. If $J_2\ne\varnothing$, then by definition of $f$ and later by~\eqref{claimst} we obtain
\begin{equation*}
\begin{split}
|\bigcup_{v\in C}B(v)| &\ge |\bigcup_{S\subseteq \{1,\dots ,\alpha -1\} : J_2\cap S\neq\varnothing}B(x_{S\cup\{\alpha\}})|+|\bigcup_{S\subseteq \{1,\dots ,\alpha -1\} : J_1\cap S\neq\varnothing\neq J_2\cap S}B(x_S)|\\
 &= \Bigl| Q \cap \Bigl(\bigcup_{i\in J_1\cup \{\alpha\}} N_i\Bigr)\cap \Bigl(\bigcup_{j\in J_2} N_j\Bigr)\Bigr| \\
&\ge |M \cap\bar N_\alpha|,
\end{split}
\end{equation*}
as desired. So we can assume $J_2=\varnothing$, in which case Lemma~\ref{parti}\emph{(\ref{vacio})} guarantees $d_S\in C$ for every $S\subseteq \{1,\dots ,\alpha -1\}$.  But then, the definition of $f$ gives $$|\bigcup_{v\in C}B(v)|\ \ge \ |\bigcup_{S\subseteq \{1,\dots ,\alpha -1\}}B(d_S)|\ =\ |M\cap\bar N_\alpha|.$$
This completes the proof of~\eqref{fini}, and thus the proof of the theorem. 
\end{proof}

\section*{Acknowledgements}
We wish to thank two anonymous referees for their careful reading and suggestions  which greatly improved the presentation of the paper.





\end{document}